\magnification=\magstep1

\centerline {\bf ON A ``REPLICATING CHARACTER 
STRING'' MODEL}
\bigskip

\noindent Richard C.\ Bradley \hfil\break
Department of Mathematics \hfil\break
Indiana University \hfil\break
Bloomington \hfil\break
Indiana 47405 \hfil\break
USA \hfil\break

\noindent E-mail address: bradleyr@indiana.edu
\hfil\break
\bigskip

   {\bf Abstract.}  In a paper of Chaudhuri and Dasgupta
published in 2006, a certain stochastic model for
``replicating character strings'' (such as in DNA
sequences) was studied.
In their model, a random ``input'' sequence was
subjected to random mutations, insertions, and deletions,
resulting in a random ``output'' sequence.
In this note, their model will be set up in a slightly
different way, in an effort to facilitate  
further development of the theory for their model.
In their 2006 paper, Chaudhuri and Dasgupta showed that
under certain conditions, strict stationarity of the
``input'' sequence would be preserved by the ``output'' sequence, and they proved a similar ``preservation'' 
result for the property of strong mixing with 
exponential mixing rate.
In our setup, we shall in spirit slightly extend their ``preservation of stationarity'' result, and also 
prove a ``preservation'' result for the property of 
absolute regularity with summable mixing rate.
\bigskip
\bigskip


\vfill\eject

\centerline {\bf 1. Introduction}
\bigskip

   Chaudhuri and Dasgupta [6] formulated and studied a 
certain stochastic model for ``replicating character strings''.
In that paper, they cited numerous other references 
where other related models had been studied, and 
in particular, they cited the book by Waterman [17] 
for the possible application of central
limit theory under strong mixing conditions in the
use of such models for the statistical analysis of 
data from biology (e.g.\ involving DNA sequences).
In this note, we shall contribute further results and techniques to the theory for the particular model in [6],
and suggest a way of setting up their model that
may allow slightly easier handling of certain 
technical details. 
   
   Let ${\bf N}$ (resp.\ ${\bf Z}$) 
denote the set of all positive integers (resp.\ the 
set of all integers).

   The model studied by Chaudhuri and Dasgupta [6]
can be briefly described as follows:
It starts with an ``input'' sequence 
$X := (X_k, k \in {\bf N})$ of random variables taking 
their values in some finite ``alphabet'' --- for 
example the set $\{A,C,G,T\}$ of letters that represent
the nucleotides in a DNA sequence.
There is another sequence 
$Z := (Z_k, k \in {\bf N})$ of random variables taking
their values in the set $\{M, I, D\}$ --- to indicate that
at a given ``time'' (or ``location'') $k$, 
there should be a ``mutation'' (M),
``insertion'' (I), or ``deletion'' (D).
(This sequence $Z$ is informally referred to below as 
the ``MID-sequence''.)\ \ 
Probabilities are assigned for what letter of the
``alphabet'' is inserted when an insertion occurs,
or what letter of the alphabet results from 
a mutation.  
(The --- perhaps high --- probability of ``no mutation''
is formally represented in this scheme as the 
probability of ``replacing a letter by itself'' when 
a ``mutation'' occurs.)\ \ 
At the end, the result is an ``output'' sequence
$Y := (Y_k, k \in {\bf N})$ of random variables, with 
the same ``alphabet'' (e.g.\ $\{A, C, G, T\}$) as the 
``input'' sequence $X$.

   In their paper, Chaudhuri and Dasgupta
[6, Theorems 3.1 and 3.2] 
established certain conditions 
under which certain properties of the
``input'' sequence --- specifically, 
strict stationarity, and strong mixing with
exponential mixing rate --- would be retained by
the ``output'' sequence.
Chaudhuri and Dasgupta [6] set up their
model using (``one-sided'') random sequences indexed by 
${\bf N}$, as described above.
In Section 3, we shall set up their model again,
but using (``two-sided'') random sequences indexed 
by ${\bf Z}$.
This will hopefully make it a little easier to handle
various technical details, such as keeping track of
relevant $\sigma$-fields when estimating mixing
rates.

   In the statements of their main results (though not
in the initial formulation of their model),
Chaudhuri and Dasgupta [6] dealt with the case where
the MID-sequence is an (irreducible, aperiodic) Markov
chain that is independent of the the ``input'' sequence.
We shall retain that ``independence'' assumption, but 
allow the MID-sequence itself to satisfy a somewhat more 
flexible dependence assumption than a ``Markov'' 
property.

   Instead of studying the mixing rates for the
``input'' and ``output'' sequences for the strong
mixing condition, we shall do so for the absolute 
regularity condition, which is stronger than strong 
mixing. 
That will provide an opportunity to illustrate the 
use of a particularly handy ``coupling'' property 
(due to Berbee [1]) that is
possessed by the absolute regularity condition but 
not by the strong mixing condition.
However, along the way, we shall also give information 
that may be relevant to the further development of the
theory for this model under the strong mixing 
condition.

   Instead of studying the case of exponential 
mixing rates (for strong mixing) as in [6], 
we shall focus on a certain
slower (``summable'') mixing rate (for absolute 
regularity) that is natural in central limit theory 
for bounded random variables (under either strong
mixing or absolute regularity). 

   In the model in [6], one
somewhat tricky facet of keeping track of relevant
$\sigma$-fields was keeping track of the changes in the
``clock'' resulting from ``deletions''.
We shall adopt an alternative technical procedure
--- switching to a new probability measure based on
conditioning on a certain event --- in the hope of
slightly simplifying that task.

   In their model, Chaudhuri and Dasgupta [6]
assumed a finite state space (for the ``input'' and
``output'' sequences), as described above.
That is the case of primary interest; but
for convenience, we shall relax that
assumption and allow the ``input'' and
``output'' sequences to consist in essence of 
real-valued random variables.
We shall actually treat those random variables as 
taking their values in $(0, \infty)$ (think of
``coding'' a real number $x$ by the positive number $e^x$),
and in an intermediate stage reserve the value 0 as a 
temporary ``place holder'' where an ``insertion'' will 
ultimately occur.

   The model in [6] directly
involved probability mass functions for what happens
when a ``mutation'' or ``insertion'' occurs.            
As a measure-theoretic convenience, we shall handle that
in a slightly different way, using independent
random variables uniformly distributed on the unit 
interval as ``randomizers''.

   In making these modifications, we shall not change
the actual model studied by Chaudhuri and Dasgupta
[6] in any significant way.
The modifications here only involve how their model
is set up.
Our quest here is in part to facilitate further 
development of the theory for their model.
There is of course the practical question, not
addressed here, of to what extent inaccuracy may occur
when, say, a ``long but finite'' DNA sequence is
modeled as a ``two-sided'' random sequence. 
 
In Section 2, some preliminary information on both the
strong mixing and absolute regularity conditions
will be given.  
In Section 3, the model in [6] 
will be spelled out with the modifications in the setup described above.
Then in Section 4, the main result of this note will be
stated and proved.
\bigskip

\centerline {\bf 2. Preliminary information on two 
mixing conditions}  
\bigskip

   In the development of the material in Sections 3 and
and 4 below, we shall start with a probability space,
and then switch to a new probability measure (on the
same measurable space) obtained
by conditioning on a certain key event.
Accordingly, in the notations in definitions below,  
the relevant probability measure will be specified
explicitly.
If only one probability measure, say $P$, is specified,
then the notation $E(\dots)$ will be tacitly understood
to mean the expected value with respect to that
particular probability measure $P$.

   Suppose $(\Omega, {\cal F})$ is a measurable space. 
Suppose $W := (W_i, i \in I)$ is a random variable/vector or stochastic process indexed by a nonempty set $I$ ---
that is, $W : \Omega \to {\bf R}^I$ is a function which is
measurable with respect to the $\sigma$-field 
${\cal F}$ on $\Omega$ and the Borel $\sigma$-field
on ${\bf R}^I$.
The $\sigma$-field ($\subset {\cal F}$) of subsets of 
$\Omega$ generated by $W$ will be denoted $\sigma(W)$ or 
$\sigma (W_i, i \in I)$.
\bigskip

   {\bf Definition 2.1.}  Suppose $(\Omega, {\cal F})$
is a measurable space, and $P$ is a probability
measure on $(\Omega, {\cal F})$.   
   
   For any two $\sigma$-fields
${\cal A}$ and ${\cal B} \subset {\cal F}$, 
define the following two measures of dependence:
$$ \eqalignno{
\alpha({\cal A}, {\cal B}; P) &:= 
\sup_{A \in {\cal A}, B \in {\cal B}}
|P(A \cap B) - P(A)P(B)|; \quad {\rm and} & (2.1) \cr 
\beta({\cal A}, {\cal B}; P) &:=
\sup {1 \over 2}\sum_{i=1}^I \sum_{j=1}^J
|P(A_i \cap B_j) - P(A_i)P(B_j)| & (2.2) \cr
} $$
where in (2.2) the supremum is taken over all pairs of
partitions $\{A_1, A_2, \dots, A_I\}$ and
$\{B_1, B_2, \dots, B_J\}$ of $\Omega$ such that
$A_i \in {\cal A}$ for each $i$ and
$B_j \in {\cal B}$ for each $j$.
It is easy to see that for any two 
$\sigma$-fields ${\cal A}$ and ${\cal B}$, one has that
$$ \alpha({\cal A}, {\cal B}; P)
\leq \beta({\cal A}, {\cal B}; P). \eqno (2.3) $$
  
   Suppose $X := (X_k, k \in {\bf Z})$ is, 
with respect to $P$, 
a strictly stationary sequence of random variables.
For each $n \in {\bf N}$ define the dependence 
coefficients
$$ \eqalignno{
\alpha(X,n;P) &:= 
\alpha(\sigma(X_k, k \leq 0), \sigma(X_k, k \geq n); P);
\ {\rm and} & (2.4) \cr
\beta(X,n;P) &:= 
\beta(\sigma(X_k, k \leq 0), \sigma(X_k, k \geq n); P). 
& (2.5) \cr 
}$$
   One trivially has that each of the sequences of numbers
$(\alpha(X,n;P),\ n \in {\bf N})$ and
\break
$(\beta(X,n;P),\ n \in {\bf N})$
is nonincreasing.
Also, by (2.3), $\alpha(X,n;P) \leq \beta(X,n;P)$
for every positive integer $n$.
The sequence $X$ is (with respect to the probability
measure $P$) ``strongly mixing'' [14] if 
$\alpha(X,n;P) \to 0$ as $n \to \infty$,
and ``absolute regular'' [16] if
$\beta(X,n;P) \to 0$ as $n \to \infty$ 
By (2.3), absolute regularity implies strong mixing.

   To motivate the results later in this paper, we shall
state a classic theorem of Ibragimov, from
Ibragimov and Linnik [10, Theorem 18.5.4].
\bigskip 

   {\bf Theorem 2.1} (Ibragimov). {\sl  
Suppose that on a probability space $(\Omega, {\cal F}, P)$,
$X := (X_k, k \in {\bf Z})$ is a strictly stationary
sequence of bounded, centered random variables
such that 
$\sum_{n=1}^\infty \alpha(X,n;P) < \infty$.

Then 
$\sigma^2 := EX_0^2 
+2 \cdot \sum_{n=1}^\infty EX_0X_n$ exists in 
$[0, \infty)$, with this sum being absolutely 
convergent. 
If further $\sigma^2 > 0$, then
$(X_1 + X_2 + \dots + X_n)/(n^{1/2} \sigma)$
converges in distribution to the $N(0,1)$ law
as $n \to \infty$.}
\bigskip

   That theorem will not be used anywhere in what follows,
but it will provide the motivation for the mathematical
development in this paper.
For example, in a statistical analysis of DNA data, 
one might deal with indicator functions 
($\{0,1\}$-valued random variables)
marking the locations of a particular
pattern of nucleotides along a DNA sequence.
Thus if strong mixing is assumed as part of the statistical 
model, then it might be natural to apply a central
limit theorem for {\it bounded\/} strongly mixing sequences
of random variables, such as Theorem 2.1.
Now the (summable) mixing rate in Theorem 2.1 is 
practically sharp. 
That was shown to be true even under absolute regularity,
by counterexamples of Davydov [7, Example 2] and 
the author [4].
The counterexample in the latter paper is a strictly
stationary, 3-state sequence that satisfies absolute
regularity with (``not quite summable'') mixing rate 
$\beta(X,n;P) = O(1/n)$.
Theorem 2.1 seems to be a natural one to use when either
strong mixing or absolute regularity is assumed in the modeling of DNA sequences;
and it is the summable mixing rate in that theorem 
that we shall focus on in this note.

   It is worth noting that Merlev\`ede and 
Peligrad [11] proved (as a special case of the main 
result in their paper) a modified, refined version
of Theorem 2.1 with the barely slower mixing rate
$\alpha(X,n;P) = o(1/n)$ and an explicit extra assumption
on the rate of growth of the variances of partial sums.

   As was mentioned in Section 1, instead of dealing 
with the strong mixing condition, 
we shall deal with absolute regularity.
That will provide an opportunity to illustrate the use
--- in Steps 5 and 6 of the proof of Lemma 4.4 in 
Section 4 ---
of a handy ``coupling'' property (from [1]) of the 
absolute regularity condition.
That property does not exist, at least in as strong a 
form, under just strong mixing.
The next three lemmas will facilitate that particular 
application of that coupling property.
\bigskip

\noindent {\bf Lemma 2.1.} {\sl
Suppose $(\Omega, {\cal F}, P)$ is a probability space,
$N$ is a positive integer, 
${\cal A}_n$ and ${\cal B}_n$, $n \in \{1, 2, \dots, N\}$
are $\sigma$-fields $\subset {\cal F}$,
and the $\sigma$-fields 
${\cal A}_n \vee {\cal B}_n$, $n \in \{1, 2, \dots, N\}$
are independent (under $P$).
Then
$$ \beta\Bigl(\bigvee_{n=1}^N {\cal A}_n,
 \bigvee_{n=1}^N {\cal B}_n ;P\Bigl) \leq
\sum_{n=1}^N \beta({\cal A}_n, {\cal B}_n;P). $$}

   In one form or another, this has long been 
part of the folklore; see e.g.\ Pinsker [12, p.\ 73].
One reference for the particular formulation here is
[5, v1, Theorem 6.2].  
(That reference also gives the exactly analogous
inequality for the dependence coefficient 
$\alpha(.\, ,.)$.)

   The next lemma has long been part of the folklore; 
but a reference for it seems hard to find.
In this lemma, the random variables $X$ and $Y$
are not assumed to be identically distributed.
In this lemma, the term ``Borel space'' means a 
measurable space $(S, {\cal S})$ that is bimeasurably
isomorphic to the space $({\bf R}, {\cal R})$ where
${\cal R}$ denotes the Borel $\sigma$-field on 
${\bf R}$. 
It is well known that ${\bf R}^{\bf N}$ (or
${\bf R}^{\bf Z}$), accompanied by its Borel
$\sigma$-field, is a Borel space.
\bigskip

\noindent {\bf Lemma 2.2.} {\sl
Suppose $(S,{\cal S})$ is a Borel space. 
Suppose $(\Omega, {\cal F}, P)$ is a probability space,
$X$ and $Y$ are random variables on this space which
take their values in $(S, {\cal S})$, and
${\cal A}$ is a $\sigma$-field $\subset {\cal F}$.
Then
$$ |\beta({\cal A}, \sigma(X); P)
- \beta({\cal A}, \sigma(Y); P)| 
\leq 2 \cdot P(X \neq Y).  \eqno (2.6) $$}

   {\it Proof.}
By symmetry, it suffices to prove that
$\beta({\cal A}, \sigma(X); P) \leq 
\beta({\cal A}, \sigma(Y); P) + 2 \cdot P(X \neq Y)$. 
Suppose $\{A_1, A_2, \dots, A_I\}$ and 
$\{B_1, B_2, \dots, B_J\}$ are each a partition of
$\Omega$, with 
$A_i \in {\cal A}$ for each $i$ and
$B_j \in \sigma(X)$ for each $j$.
It suffices to show that
$$ {1 \over 2} \sum_{i=1}^I \sum_{j=1}^J
|P(A_i \cap B_j) - P(A_i)P(B_j)| \leq
\beta({\cal A}, \sigma(Y); P) + 2P(X \neq Y). 
\eqno (2.7) $$

By a well known measure-theoretic fact (a standard
elementary generalization of [2, Theorem 20.1]),
there exists a partition $\{S_1, S_2, \dots, S_J\}$
of $S$ with $S_j \in {\cal S}$ for each $j$, 
such that $B_j = \{X \in S_j\}$ for each $j$.      
    
For any event $A$, 
$$ \eqalignno{
\sum_{j=1}^J&
|P(A \cap \{X \in S_j\}) - P(A \cap \{Y \in S_j\})| \cr
&\leq 
\sum_{j=1}^J
|P(A \cap \{X \in S_j\} \cap \{X = Y\})
 +P(A \cap \{X \in S_j\} \cap \{X \neq Y\}) \cr
 & \indent\indent -P(A \cap \{Y \in S_j\} \cap \{X = Y\})
 -P(A \cap \{Y \in S_j\} \cap \{X \neq Y\})| \cr
&=
\sum_{j=1}^J
|P(A \cap \{X \in S_j\} \cap \{X \neq Y\}) 
 -P(A \cap \{Y \in S_j\} \cap \{X \neq Y\})| \cr
 & \leq 2P(A \cap \{X \neq Y\}). \cr  
}$$
Applying that with $A = A_i$ and then also with 
$A = \Omega$, one has that
$$ \eqalignno{
\sum_{i=1}^I & \sum_{j=1}^J
|P(A_i \cap B_j) - P(A_i)P(B_j)| \cr
& \leq \sum_{i=1}^I \sum_{j=1}^J
|P(A_i \cap \{X \in S_j\}) - P(A_i \cap \{Y \in S_j\})| \cr 
& \indent + \sum_{i=1}^I \sum_{j=1}^J
|P(A_i \cap \{Y \in S_j\}) - P(A_i)P(Y \in S_j)| \cr
& \indent + \sum_{i=1}^I \sum_{j=1}^J
|P(A_i)P(Y \in S_j) - P(A_i)P(X \in S_j)| \cr
& \leq \Bigl[\sum_{i=1}^I 2P(A_i \cap \{X \neq Y\}) \Bigl] 
+ 2\beta({\cal A}, \sigma(Y); P) \cr
& \indent + \sum_{i=1}^I \Bigl[P(A_i) \sum_{j=1}^J 
|P(Y \in S_j) - P(X \in S_j)|\Bigl]\cr
& \leq 2P(X \neq Y) + 2\beta({\cal A}, \sigma(Y); P) 
+ \sum_{i=1}^I [P(A_i) \cdot 2P(X \neq Y)] \cr 
&=  4P(X \neq Y) + 2\beta({\cal A}, \sigma(Y); P). \cr   
}$$
Thus (2.7) holds. 
That completes the proof.
\bigskip

   The following lemma, the final item of Section 2
here, will play a key role (in Section 4) in the 
comparison of the mixing rates for the ``input'' and ``output'' sequences in the model in [6]. 
It will be applied for absolute regularity; but it holds
under strong mixing (as stated and proved here)
as well.
\bigskip 

\noindent {\bf Lemma 2.3.}  {\sl
Suppose $(H(0), H(1), H(2), H(3), \dots)$ is a nonincreasing
sequence of nonnegative numbers such that
$\sum_{n=0}^\infty H(n) < \infty$.
Suppose that on some probability space 
$(\Omega, {\cal F}, P)$,
$X := (X_k, k \in {\bf Z})$ is a nondegenerate,
strictly stationary sequence of random variables
taking only the values 0 and 1, such that
$\sum_{n=1}^\infty \alpha(X,n;P) < \infty$.
Then
$$ \sum_{n=1}^\infty EH(X_1 + X_2 + \dots + X_n) < \infty. 
\eqno (2.8) $$}

   {\it Proof.} 
Referring to the hypothesis (of Lemma 2.3), define the
number 
$$   p := P(X_0 = 1) = EX_0 > 0  \eqno (2.9) $$

   Define the constant random variable $S_0 := 0$; and
for each positive integer $n$, define the partial sum
$S_n := X_1 + X_2 + \dots + X_n$.
By the hypothesis (of Lemma 2.3), the sequence $X$ is
strongly mixing and hence ergodic.
Hence from (2.9), one has that  
$S_n \to \infty$ (monotonically) a.s.\ as $n \to \infty$.
For technical convenience, without loss of generality (redefining the random variables $X_k$ on a $P$-null set 
if necessary), we assume that that happens at literally 
every $\omega \in \Omega$.

   For each nonnegative integer $j$, define the random 
variable 
$$\eta_j := {\rm card}\{n \in {\bf N}: S_n = j\}.   
\eqno (2.10) $$
Then for every integer $J \geq 0$,
$$ \sum_{j=0}^J \eta_j = \max \{n \geq 0: S_n = J\}.
\eqno (2.11) $$

   In what follows, for any real number $x$, let
$[x]$ denote the greatest integer $\leq x$.
Also, in the calculations below, by the hypothesis 
(of Lemma 2.3), all sums
and summands (``numerical'' or random) 
take their values in
$[0,\infty] := [0,\infty) \cup \{\infty\}$, 
and hence one
can change the orders of summations arbitrarily.

   Recall that for any nonnegative integer-valued
random variable $W$, 
$EW = \sum_{n=1}^\infty P(W \allowbreak \geq n)$.
For each integer $J \geq 0$, by (2.11) and
the trivial inequality $P(S_n \leq J) \leq 1$,
one has that
$$
E\Bigl(\sum_{j=0}^J \eta_j\Bigl) 
= \sum_{n=1}^\infty P\Bigl(\sum_{j=0}^J \eta_j \geq n\Bigl)
= \sum_{n=1}^\infty P(S_n \leq J)
\leq 2J/p + \sum_{n = [2J/p]+1}^\infty P(S_n \leq J).
\eqno (2.12)
$$ 

   Let us examine the last sum in (2.12).
For each integer $J \geq 0$ and each integer $n > 2J/p$, 
one has that $J < np/2$ and hence $J-np < -np/2$. 
Hence for each integer $J \geq 0$,
$$ \eqalignno{
\sum_{n=[2J/p]+1}^\infty &P(S_n \leq J)
\leq \sum_{n=[2J/p]+1}^\infty P(S_n -np \leq -np/2)
\leq \sum_{n=[2J/p]+1}^\infty P(|S_n -np| \geq np/2) \cr
& \leq \sum_{n=[2J/p]+1}^\infty
(np/2)^{-4} E(S_n -np)^4 
\leq  (16/p^4) \sum_{n=1}^\infty n^{-4}E(S_n-np)^4. 
& (2.13) \cr 
}$$

   Extend the definition (2.4) to include $n = 0$ there.
Then here for each positive integer $n$, 
$$ E(S_n - np)^4 \leq 
(20,000)n \cdot \sum_{m=0}^{n-1} (m+1)^2 \alpha(X,m;P)
+ 24n^2 \Bigl[\sum_{m=0}^{n-1} \alpha(X,m;P)\Bigl]^2. 
\eqno (2.14)$$
This is a simple direct application of 
[5, v2, Theorem 14.63],
which in turn is a convenient but crude version of a much sharper and more general inequality due to 
Rio [13, Th\'eor\`eme 2.1].
Keep in mind that since the random variables $X_k$ take only the values 0 and 1, the (``upper-tail'') quantile functions
in those particular statements in both references take only the values 0 and 1.
Eq.\ (2.14) can also be obtained directly, in a 
sharper form, from a careful
calculation from Ibragimov's proof in 
[10, Theorem 18.5.4] of Theorem 2.1  
(examine carefully the argument for 
[10, Lemma 18.5.2]).

   Our next task is to use (2.14) to show that the last
sum in (2.13) is finite.
From simple calculus, let $C_1$ be a positive 
number such that
$\sum_{n=q}^\infty n^{-3} \leq C_1q^{-2}$ 
for every positive integer $q$. 
By the hypothesis (of Lemma 2.3),
$$ \eqalignno{
\sum_{n=1}^\infty \Bigl[ n^{-4} \cdot n \cdot 
\sum_{m=0}^{n-1} (m+1)^2 \alpha(X,m;P) \Bigl] 
& = \sum_{m=0}^\infty \sum_{n=m+1}^\infty 
n^{-3} (m+1)^2 \alpha(X,m;P) \cr
& \leq \sum_{m=0}^\infty C_1 \alpha(X,m;P) 
< \infty. \cr 
}$$
Also, trivially by the hypothesis,
$ \sum_{n=1}^\infty \Bigl[n^{-4} \cdot
n^2 \Bigl[\sum_{m=0}^{n-1} \alpha(X,m;P)\Bigl]^2 \Bigl]
< \infty$. 
Applying those two inequalities to (2.14), one obtains
that the last sum in (2.13) is finite.

   Accordingly, defining the finite numbers
$C_2 := (16/p^4) \sum_{n=1}^\infty n^{-4}E(S_n-np)^4$
and $C_3 := (2/p) + C_2$,
one has by (2.12) and (2.13) that for every integer 
$J \geq 0$,
$$ E\Bigl(\sum_{j=0}^J \eta_j\Bigl) \leq 2J/p + C_2
\leq C_3(J+1). \eqno (2.15) $$

   Now refer to the function $H$ in the statement
of Lemma 2.3.
By the hypothesis, $H(n) \downarrow 0$ as $n \to \infty$. 
Using the notation $S(n)$
for $S_n$ in subscripts, one has the equality of
nonnegative random variables 
(possibly taking the value $\infty$)
$$ \eqalignno{ 
\sum_{n=1}^\infty H(S_n)
&= \sum_{j=0}^\infty \sum_{\{n \in {\bf N}: S(n) = j\}} H(j)
= \sum_{j=0}^\infty H(j) \cdot \eta_j \cr
& = \sum_{j=0}^\infty 
\sum_{i=j}^{\infty} \eta_j \cdot [H(i)-H(i+1)] 
= \sum_{i=0}^\infty \sum_{j=0}^i [H(i)-H(i+1)] 
\cdot \eta_j.    
}$$
Hence by (2.15),
$$ \eqalignno{ 
E \sum_{n=1}^\infty H(S_n)
&\leq \sum_{i=0}^\infty \Bigl[\, [H(i)-H(i+1)] 
\cdot C_3(i+1)\Bigl] 
= C_3 \sum_{i=0}^\infty \sum_{j=0}^i [H(i)-H(i+1)] \cr
&= C_3 \sum_{j=0}^\infty \sum_{i=j}^\infty [H(i)-H(i+1)]
= C_3 \sum_{j=0}^\infty H(j) < \infty. \cr    
\cr 
}$$
Thus (2.8) holds.  
That completes the proof of Lemma 2.3. 
\bigskip

\centerline {\bf 3. The model of Chaudhuri and Dasgupta, 
in ``two-sided'' form}
\bigskip

   In this section, we shall spell out, step by step,
the ``replicating character string'' model studied by
Chaudhuri and Dasgupta [6].
As explained in Section 1, essentially the only 
changes here will be in the 
``style'': 
(i) the use of (``two-sided'') random sequences indexed 
by ${\bf Z}$, 
rather than (``one-sided'') random sequences 
indexed by ${\bf N}$, and 
(ii) the trivial allowing of the ``alphabet'' or ``state space'' to be $(0,\infty)$ instead of just a finite set.
In the presentation here, the essential mathematical 
substance of their model will not be changed at all.
Much of the notations below will be taken directly 
from their paper.
For convenient reference, the stages in this
construction will be referred to as 
paragraphs (P1), (P2), etc. 
\bigskip

   {\bf (P1).} Suppose $(\Omega, {\cal F}, P)$ is a probability space.
All random variables defined below, will be understood
to be defined on this space.
\bigskip
 
   {\bf (P2).} Suppose $X := (X_k, k \in {\bf Z})$ is 
(under $P$) a strictly stationary sequence of random 
variables taking their values in the open half line $(0,\infty)$.

   (This is the ``input'' sequence, as in the model 
in [6].)
\bigskip

   {\bf (P3).} Suppose $Z := (Z_k, k \in {\bf Z})$ is 
(under $P$) a strictly stationary, ergodic sequence of 
random variables taking their values in the set 
$\{M,I,D\}$, with this sequence $Z$ being 
independent of the sequence $X$.
Assume further that $P(Z_0 = s) > 0$ for all three 
elements $s \in \{M,I,D\}$.
For technical convenience, without loss of generality
(by ergodicity), 
assume that for {\it every\/} $\omega \in \Omega$ 
and all three elements $s \in \{M,I,D\}$,
$Z_k(\omega) = s$ for infinitely many negative
integers $k$ and infinitely many positive integers $k$. 
 
   (Again, the letters $M$, $I$, and $D$, stand
for ``mutation'', ``insertion'', and ``deletion'';
$Z$ is the ``MID-sequence'', as in [6].)
\bigskip

   {\bf (P4).} Define the strictly increasing sequence
$\zeta := (\zeta_j, j \in {\bf Z})$ 
of integer-valued random variables as follows:
For every $\omega \in \Omega$,
$$ \eqalignno{ 
\dots &< \zeta_{-2}(\omega) < \zeta_{-1}(\omega)
< \zeta_0(\omega) \leq 0 < 1 \leq 
\zeta_1(\omega) < \zeta_2(\omega) < \zeta_3(\omega)
< \dots \quad {\rm and} \cr
& \{j \in {\bf Z}: Z_j(\omega) \in \{M,D\} \}
= \{ \zeta_k(\omega): k \in {\bf Z} \}. & (3.1) \cr
} $$
The random variables $\zeta_k$ will
sometimes be written $\zeta(k)$ for typographical
convenience.
\bigskip

   {\bf (P5).} Define the sequence 
$\overline X := (\overline X_k, k \in {\bf Z})$
of random variables 
(taking their values in the closed half line $[0,\infty)$)
as follows:
For every $\omega \in \Omega$,
$$ 
\forall k \in {\bf Z},\  
\overline X_{\zeta(k)(\omega)}(\omega) := X_k(\omega);
\quad {\rm and} \quad 
\forall j \notin \{\zeta_k(\omega): k \in {\bf Z}\},\ 
\overline X_j(\omega) := 0.  \eqno (3.2)    
$$

   The state 0 is used here only as a ``temporary
placeholder'' for a spot where an ``insertion'' will
eventually occur (in eq.\ (3.5) below, in the case
where $\overline Y_\ell(\omega) = 0$ there).
That was the sole motivation for choosing for the
original sequence $X$ a state space, namely $(0, \infty)$,
that does not include 0.
\bigskip

   {\bf (P6).} Define the strictly increasing sequence
$\xi := (\xi_j, j \in {\bf Z})$ 
of integer-valued random variables as follows:
For every $\omega \in \Omega$,
$$ \eqalignno{ 
\dots &< \xi_{-2}(\omega) < \xi_{-1}(\omega)
< \xi_0(\omega) \leq 0 < 1 \leq 
\xi_1(\omega) < \xi_2(\omega) < \xi_3(\omega)
< \dots \quad {\rm and} \cr
& \{j \in {\bf Z}: Z_j(\omega) \in \{M,I\} \}
= \{ \xi_k(\omega): k \in {\bf Z} \}. & (3.3) \cr
} $$
These random variables $\xi_k$ will
sometimes be written $\xi(k)$. 
\bigskip

   {\bf (P7).} Define the sequence 
$\overline Y := (\overline Y_\ell, \ell \in {\bf Z})$
of random variables 
(taking their values in $[0,\infty)$)
as follows:
For every $\omega \in \Omega$,
$$ 
\forall \ell \in {\bf Z},\
\overline Y_\ell(\omega) := 
\overline X_{\xi(\ell)(\omega)}(\omega).\eqno (3.4)    
$$

   {\bf (P8).} (i) Let $U := (U_k, k \in {\bf Z})$ be 
(under $P$) a sequence of independent, identically 
distributed random variables, each uniformly distributed 
on the interval $[0,1]$, with this sequence $U$ being independent of the pair of sequences $(X,Z)$.

   (ii) Let $g: (0,\infty) \times [0,1] \to (0,\infty)$ 
be a Borel function.

   (iii) Let $h: [0,1] \to (0,\infty)$ be a Borel 
function.

   (iv) Define the sequence 
$Y := (Y_\ell, \ell \in{\bf Z})$
of random variables 
(taking their values in $(0,\infty)$)
as follows:  
For every $\omega \in \Omega$,
$$ \forall \ell \in {\bf Z},\  
Y_\ell(\omega) := \cases{
g(\overline Y_{\ell}(\omega), U_\ell(\omega)) & if 
$\overline Y_{\ell}(\omega) \in (0,\infty)$ \cr
h(U_\ell(\omega)) & if $\overline Y_{\ell}(\omega) = 0$. \cr
} \eqno (3.5) $$

   (This sequence $Y$ is the ``output'' sequence, as 
in [6].)

   The final two ``paragraphs'' below give a few more
items that were not needed in the formulation of the
``output'' sequence $Y$ but will be needed in the
formulation of the main result (Theorem 4.1 
in Section 4). 
\bigskip

   {\bf (P9).}  Let $P_0$ denote the probability measure
on $(\Omega, {\cal F})$ defined as follows:
$$  \forall F \in {\cal F},\  
P_0(F) := P(F|Z_0 \in \{M,I\}) = P(F|\xi_0 = 0). 
\eqno (3.6) $$
(The second equality follows from paragraph (P6).)
\bigskip

   {\bf (P10).} Define the sequence 
$V := (V_k, k \in {\bf Z})$ of 
$(\{M,I\} \times {\bf N})$-valued 
random variables as follows:
$$ \forall k \in {\bf Z},\ 
V_k := (Z_{\xi(k)}, \xi_k - \xi_{k-1}). \eqno (3.7) $$
Also, define the sequence 
$\Upsilon := (\Upsilon_k, k \in {\bf Z})$ of
$(\{M,I\} \times {\bf N} \times (0,\infty))$-valued 
random variables as follows:
$$ \forall k \in {\bf Z},\ 
\Upsilon_k := (V_k, Y_k) = 
(Z_{\xi(k)}, \xi_k - \xi_{k-1}, Y_k). \eqno (3.8) $$ 
That completes the (``two-sided'') presentation of the 
model in [6]. \bigskip  

\noindent {\bf Remark 3.1.} Here are several comments
pertaining to the model from [6] as spelled out in
paragraphs (P1)-(P10) above.
\smallskip

   (A)  It is well known from renewal theory that even 
though the sequence $Z$ is (under the original
probability measure $P$) strictly stationary, 
the sequence $(Z_{\xi(k)}, k \in {\bf Z})$ is in 
general {\it not\/} strictly stationary under $P$.
As a consequence, under $P$ the ``output'' sequence $Y$
will in general not be strictly stationary.
To obtain the stationarity of $Y$, 
Chaudhuri and Dasgupta [6, Theorem 3.1] 
assumed directly that that sequence 
$(Z_{\xi(k)}, k \in {\bf Z})$
(though not necessarily the entire sequence $Z$) 
is strictly stationary, 
with the original MID-sequence $Z$ itself being a 
Markov chain with certain properties.
Theorem 4.1 in Section 4 (the main result of this note)
will employ the procedure, 
common in renewal theory, of formally 
switching to the new probability measure $P_0$ in (3.6),
which under our assumptions will yield the 
stationarity of $Y$ (and of the
entire sequence $\Upsilon$ in (3.8)), with no
``Markov'' assumption.
That is the role here of the probability
measure $P_0$. 
\smallskip

   (B) By paragraphs (P4), (P6), and (P10), one has
that $\sigma(\zeta, \xi, V) \subset \sigma(Z)$. 
Recall from paragraphs (P3) and (P8) that 
under $P$, the random sequences $X$, $Z$, and $U$ are
independent of each other.
By (3.6) and a trivial argument, those three
sequences are independent of each other under $P_0$
as well.
Also by (3.6), the random sequences $X$ and $U$ 
(but in general not $Z$ or even $V$) each have the same 
distribution under $P_0$ as they do under $P$.
In particular,
$$ \forall n \in {\bf N},\ 
\beta(X,n;P_0) = \beta(X,n;P) \eqno (3.9) $$
(and the analogous equality holds for $\alpha(.\, ,.)$).
\smallskip

   (C) In paragraph (P3), it was implicitly understood
in the phrase ``without loss of generality'' 
that on a certain ``bad'' event $F$ with $P(F) = 0$,
one might need to redefine certain random variables 
$Z_k, k \in {\bf Z}$.
By (3.6) $P_0(F) = 0$ as well, and hence that phrase
``without loss of generality'' applies under $P_0$
as well as under $P$.
\smallskip

   (D) In [6] (with finite ``alphabet''), 
the probabilities involving ``mutations'' and 
``insertions'' were specified directly.
In [6] it was also pointed out how the context of
``mutation'' could include, as part of the model, 
``high probability of no mutation''. 
Paragraph (P8) just gives an alternative 
way to set all that up, using the independent 
random variables $U_k$ uniformly distributed on the 
interval $[0,1]$ as ``randomizers'', and using 
appropriate choices of the Borel functions $g$ 
(to determine ``mutations'') and $h$
(to determine ``insertions''). 
In particular, for a given $x \in (0,\infty)$, 
the function $g(x,u), u \in [0,1]$ can be defined to be
equal to $x$ itself (``no mutation'') for $u$ in ``most''
of the interval $[0,1]$.     
\bigskip

\centerline {\bf 4. The main result and its proof}
\bigskip

   This section is devoted to the proof of the following
theorem, the main result of this note:
\bigskip  

\noindent {\bf Theorem 4.1.} {\sl 
Assume the entire context of paragraphs (P1)-(P10), 
with all assumptions there satisfied.
  
   (I) Under the probability measure $P_0$ in (3.6),
the sequence $\Upsilon$ (in paragraph (P10)) 
is strictly stationary (and hence under $P_0$ 
the sequences $V$ and $Y$ are each strictly stationary).

   (II) If also $\sum_{n=1}^\infty \beta(X,n;P) < \infty$
(see also (3.9)) and 
$\sum_{n=1}^\infty \beta(V,n;P_0) \allowbreak < \infty$, 
then 
$\sum_{n=1}^\infty \beta(Y,n;P_0) \leq 
\sum_{n=1}^\infty \beta(\Upsilon,n;P_0) <\infty$.}
\bigskip

   Statement (I) is in spirit a slight extension of
[6, Theorem 3.1], which in their setup was a 
corresponding ``preservation of strict stationarity''
result involving a ``Markov'' assumption on the 
MID-sequence. 
Statement (II) was inspired by [6, Theorem 3.2],
which in their setup was a corresponding 
``preservation of mixing rate'' result involving 
strong mixing with exponential mixing rate.
It seems clear that the setup here in 
paragraphs (P1)-(P10), involving ``two-sided'' random sequences, can facilitate the proofs of such 
``preservation of mixing rates'' results involving 
absolute regularity, such as Statement (II) here;
but it is yet to be determined to what extent the 
setup here might facilitate such 
results involving strong mixing.         
\bigskip

   We shall first prove Statement (I).
The proof given below will be a somewhat modified 
version of the argument for [6, Theorem 3.1].
The argument will proceed through a series of lemmas.
The first lemma is of a standard form.
(In closely related contexts, a very similar fact 
was used in [3, proof of Lemma 5] and in [4, pp.\ 7-8];
see also [5, v3, Theorem 26.4(I)].)
\bigskip

\noindent {\bf Lemma 4.1.} {\sl
In the context of paragraphs (P1)-(P10) (with all
assumptions there satisfied), the sequence
$((Z_k, \overline X_k), k \in {\bf Z})$ is, under the
probability measure $P$, strictly stationary.}
\bigskip

   {\it Sketch of proof.}
Suppose $j$ is any integer.  
Define the integer-valued random variable
$T := \max\{k \in {\bf Z}: \zeta_k \leq j\}$.
The entire array 
$((Z_k, k \geq j+1), (\overline X_k, k \geq j+1))$
can be represented as
$\phi((Z_k, k \geq j+1),
(X_{T+1}, X_{T+2}, X_{T+3}, \dots))$,
where the (measurable) function 
$\phi : \{M,I,D\}^{\bf N} \times (0,\infty)^{\bf N}
\to \{M,I,D\}^{\bf N} \times [0,\infty)^{\bf N}$ 
does not depend on $j$.
Under $P$, regardless of $j$, by the assumptions in 
paragraphs (P2)-(P5) and an elementary argument,
the sequence 
$(X_{T+1}, \allowbreak X_{T+2}, X_{T+3}, \dots)$
is independent of $\sigma(Z,T)$ ($= \sigma(Z)$) 
and has the same distribution as the sequence 
$(X_1, X_2, X_3, \dots)$,
and Lemma 4.1 then follows easily. 
\bigskip

\noindent {\bf Lemma 4.2.} {\sl 
Suppose $L \geq 3$ is an integer.
Suppose that for each $\ell \in \{1,2,\dots, L\}$,
$s_\ell \in \{M,I\}$,
$N_\ell$ is a positive integer, and
$B_\ell$ is a Borel subset of $[0,\infty)$.
For each $J \in \{1,2,\dots,L-1\}$, define the event
(see (3.7) and (3.3)) 
$$ F_J := \Bigl\{Z_0 \in \{M,I\}\Bigl\} \bigcap
\Bigl[\bigcap_{\ell = 1}^L \Bigl(\{V_{-J+\ell} = 
(s_\ell, N_\ell)\}
\cap \{\overline Y_{-J+\ell} \in B_\ell\}\Bigl)\Bigl].
\eqno (4.1) $$ 
Then $P(F_1) = P(F_2) = \dots = P(F_{L-1})$.}
\bigskip

   {\it Proof.} Define the integers $m_0 := 0$ and
$m_\ell := N_1 + N_2 + \dots + N_\ell$ for
$\ell \in \{1, 2, \dots, L\}$.
These integers $m_\ell$ will sometimes be written
below as $m(\ell)$.
Define the set 
$S := \{1,2,\dots, \allowbreak m_L\} - 
\{m_1,m_2, \dots, m_L\}$. 

   Suppose $J \in \{1, 2, \dots, L-1\}$.
   
   By (3.3), 
$\{Z_0 \in \{M,I\}\} = \{\xi_0 = 0\}$.
As a consequence, one has the equality of events
$$ \eqalignno{
\Bigl\{Z_0 \in &\{M,I\}\Bigl\} \bigcap 
\Bigl[\bigcap_{\ell = 1}^L 
\{\xi_{-J+\ell} - \xi_{-J+\ell-1}= N_\ell \}\Bigl]
= \bigcap_{\ell=0}^L \{\xi_{-J+\ell} = -m_J+m_\ell\} 
\cr
& = \Bigl[ \bigcap_{\ell = 0}^L 
\{Z_{-m(J) + m(\ell)} \in \{M,I\} \Bigl] \bigcap 
\Bigl[ \bigcap_{u \in S} \{Z_{-m(J)+u} = D\} \Bigl].
& (4.2)
\cr 
}$$
Referring to (4.1) and applying both equalities in (4.2) carefully, one obtains
$$\eqalignno{
F_J &=  
\Bigl\{Z_0 \in \{M,I\}\Bigl\} \bigcap \cr 
& \indent \quad \Bigl[\bigcap_{\ell = 1}^L
\Bigl( \{Z_{\xi(-J+\ell)} = s_\ell \} \cap
\{\xi_{-J+\ell} - \xi_{-J+\ell-1}= N_\ell \}
\cap \{\overline X_{\xi(-J+\ell)} \in B_\ell \} \Bigl)
\Bigl]
\cr
&= \Bigl[\bigcap_{\ell=0}^L \{\xi_{-J+\ell} = 
-m_J+m_\ell\} \Bigl] \bigcap \Bigl[\bigcap_{\ell = 1}^L
\Bigl( \{Z_{-m(J)+ m(\ell)} = s_\ell \}
\cap \{\overline X_{-m(J)+ m(\ell)} \in B_\ell\} \Bigl)
\Bigl] 
\cr
&= \Bigl\{ Z_{-m(J)} \in \{M,I\} \Bigl\} \bigcap
\biggl[ \bigcap_{\ell = 1}^L
\Bigl( \{Z_{-m(J)+ m(\ell)} = s_\ell \}
\cap \{\overline X_{-m(J)+ m(\ell)} \in B_\ell\} \Bigl)
\biggl] \cr
& \indent \quad \bigcap 
\Bigl[ \bigcap_{u \in S} \{Z_{-m(J)+u} = D\} \Bigl].
& (4.3)   
}$$ 
By Lemma 4.1, the probability (under $P$) of the last expression in (4.3) does not
depend on $J$  ($\in \{1, 2, \dots, L-1\}$). 
Thus Lemma 4.2 holds.
\bigskip

\noindent {\bf Lemma 4.3.} {\sl 
The sequence 
$((V_\ell, \overline Y_\ell), \ell \in {\bf Z})$ of
$((\{M,I\} \times {\bf N}) \times [0,\infty))$-valued
random variables is strictly stationary under $P_0$.}
\bigskip

   {\it Proof.} 
Suppose $j \in {\bf Z}$ and $n \in {\bf N}$.
It suffices to prove that under $P_0$, the 
``random vectors'' 
$((V_\ell, \overline Y_\ell), 
\ell \in \{j+1, j+2, \dots, j+n\})$
and 
$((V_\ell, \overline Y_\ell), 
\ell \in \{j+2, j+3, \dots, j+n+1\})$     
have the same distribution
(on $(\{M,I\}\times {\bf N} \times [0,\infty))^n$).

   Let $J$ and $L$ be positive integers such that
$\{J-1,J\} \subset \{1,2,\dots, L-1\}$ and
$\{j+1, j+2, \dots, j+n\} \subset 
\{-J+1, -J+2, \dots, -J+L\}$.   
It suffices to prove that under $P_0$, the ``random
vectors''
$((V_\ell, \overline Y_\ell), 
\ell \in \{-J+1, -J+2, \dots, -J+L\})$
and 
$((V_\ell, \overline Y_\ell), 
\ell \in \{-J+2, -J+3, \dots, -J+L+1\})$
have the same distribution.
But that holds by (3.6), Lemma 4.2, and a trivial 
calculation.
Thus Lemma 4.3 holds.
\bigskip  

   {\it Proof of Statement (I) in Theorem 4.1.}
By paragraph (P8) (see Remark 3.1(B)), the sequence $U$
is, under $P_0$, independent of the sequence
$((V_\ell, \overline Y_\ell), \ell \in {\bf Z})$.
It follows from paragraph (P8) and Lemma 4.3 (again
see Remark 3.1(B)) that under $P_0$,
the sequence 
$((V_\ell, \overline Y_\ell, U_\ell), \ell \in {\bf Z})$ 
is strictly stationary.
Now Statement (I) in Theorem 4.1 holds by (3.5) and (3.8).
\bigskip

   The proof of Statement (II) in Theorem 4.1 will be
based on the following lemma.
In what follows, $E_0(\dots)$ denotes expected value
with respect to the probability measure $P_0$.
The indicator function of a given event $A$ will be
denoted $I(A)$.
\bigskip

\noindent {\bf Lemma 4.4.} {\sl
In the context of (P1)-(P10) (with all assumptions
there satisfied), suppose also that
$\sum_{n=1}^\infty \beta(X,n;P) < \infty$.
Define the sequence 
$(H(n), n \in \{0\} \cap {\bf N})$ of
nonnegative numbers 
as follows:
For each $n \geq 0$, $H(n) := \beta(X,n+1;P)$.

Suppose $N$ is an integer such that $N \geq 2$.
Then
$$ \beta(\Upsilon,N;P_0) \leq
\beta(V,N;P_0) + 
2E_0 H\Bigl(\sum_{i=1}^{N-1} I(Z_{\xi(i)} = M)\Bigl). 
\eqno (4.4)$$} 

   {\it Proof.}
The proof of this lemma will proceed through a series of
``steps''.
\bigskip

   {\it Step 1.} 
Refer to the integer $N \geq 2$ in the hypothesis 
(of Lemma 4.4). 
Define the nonnegative integer-valued
random variable $T$ by
$$ \eqalignno{ 
T &:= {\rm card}\Bigl\{k \in {\bf N}: 
1 \leq k \leq \xi_N-1\
{\rm and}\ Z_k \in \{M,D\} \Bigl\} \cr
&\, = \max\{j \in \{0\}\cup{\bf N}: \zeta_j < \xi_N\}.
& (4.5) \cr
}$$
(The second equality in (4.5) holds by (3.1).)\ \ 
Define the (``one-sided'') sequence
$X^* := (X_1^*, X_2^*,X_3^*, \dots)$ of random variables
as follows:
$$ \forall k\geq 1,\ X_k^* := X_{k+T}. \eqno (4.6) $$

   {\it Step 2.}
Now let us first look at the random variable 
$\overline Y_N$.
Suppose $\omega \in \Omega$.
If $Z_{\xi(N)(\omega)}(\omega) = I$, then
$\xi_N(\omega) \notin \{\zeta_k(\omega): k \in {\bf Z}\}$
by (3.1), and 
$\overline Y_N(\omega) = 0$ by (3.4) and (3.2).
If instead $Z_{\xi(N)(\omega)}(\omega) = M$ (the only
other possibility, by (3.3)), then
for some $q \geq 1$, 
$\xi_N(\omega) = \zeta_q(\omega)$, hence
$q = T(\omega)+1$ by (4.5), and hence
$\overline Y_N(\omega) = 
\overline X_{\zeta(q)(\omega)}(\omega) = X_q(\omega)
= X_{T(\omega)+1}(\omega) = X_1^*(\omega)$ 
by (3.4), (3.2) and (4.6).

Thus
$ \overline Y_N = 
0 \cdot I(Z_{\xi(N)} = I) 
+ X_1^* \cdot I(Z_{\xi(N)} = M)$.
Hence by (3.7),
$$ \sigma(\overline Y_N) \subset \sigma (V_N, X_1^*). \eqno (4.7)$$ 

   {\it Step 3.}
Now suppose $\ell$ is any integer such that $\ell > N$.
Our task here in Step 3 is to obtain some sort of analog
of (4.7) for $\overline Y_\ell$.

   First define the random variable
$$ \tau := {\rm card}\Bigl\{k \in {\bf N}: 
\xi_N \leq k \leq \xi_\ell \
{\rm and}\ Z_k \in \{M,D\} \Bigl\}. \eqno (4.8) $$ 
Then
$\tau = [\sum_{i=N}^\ell I(Z_{\xi(i)} = M)] 
+ [\sum_{i=N+1}^\ell (\xi_i - \xi_{i-1} -1)]$.
Here in the right hand side, 
for a given $\omega \in \Omega$, by (3.3),
the first sum is simply the number of indices $k$
in the set in (4.8) such that $Z_k(\omega) = M$,
and the second sum is simply the number of 
indices $k$ in that set such that $Z_k(\omega) = D$.
(Either sum can be zero.)\ \ 
From that expression for $\tau$, one has by (3.7) that
$$ \sigma(\tau) \subset 
\sigma(V_N, V_{N+1}, \dots, V_\ell).
\eqno (4.9) $$ 

   Now suppose $\omega \in \Omega$.
Consider first the case where 
$Z_{\xi(\ell)(\omega)}(\omega) = M$.
Then for some $q \geq 1$, 
$\xi_\ell(\omega) = \zeta_q(\omega)$;
and by (4.5), (4.8), and (3.1),
$q = T(\omega) + \tau(\omega)$.
Hence by (3.4), (3.2), and (4.6), 
$\overline Y_\ell(\omega) = 
\overline X_{\zeta(q)(\omega)}(\omega) 
= X_q(\omega) = X^* _{\tau(\omega)}(\omega)$.
Also in the case where 
$Z_{\xi(\ell)(\omega)}(\omega) = M$,
one has that $\tau(\omega) \geq 1$ by (4.8).
If instead $Z_{\xi(\ell)(\omega)}(\omega) = I$
(the only other possibility), then
$\overline Y_\ell(\omega) = 0$ by (3.4) and (3.2).

   Thus (for our given $\ell > N$), putting all 
those pieces together,
$$ \overline Y_\ell = 
0 \cdot I(Z_{\xi(\ell)} = I) +
X^*_\tau \cdot I(Z_{\xi(\ell)} = M)
= 0 + \sum_{t=1}^\infty [X^*_t 
\cdot I(\tau = t) \cdot I(Z_{\xi(\ell)} = M)], 
$$
and hence by (4.9) and (3.7),
$$ \sigma(\overline Y_\ell) \subset
\sigma(V_N, V_{N+1}, \dots, V_\ell)
\vee \sigma(X^*_1, X^*_2, X^*_3, \dots).  
\eqno (4.10) $$

   {\it Step 4.} 
Combining (4.7) and (4.10), one now has that
$$ \sigma(\overline Y_N, \overline Y_{N+1}, 
\overline Y_{N+2}, \dots)
\subset
\sigma(V_N, V_{N+1}, V_{N+2}, \dots)
\vee \sigma(X^*_1, X^*_2, X^*_3, \dots).
\eqno (4.11) $$

   Define the two random arrays ${\bf A}$ and 
${\bf B}$ by
$$ \eqalignno{  
{\bf A} &:= \Bigl((X_k, k \leq 0); (V_k, k \leq 0);
              (U_k, k \leq 0)\Bigl) \quad {\rm and}\ 
              & (4.12) \cr   
{\bf B} &:= \Bigl(X^*; (V_k, k \geq N);
              (U_k, k \geq N)\Bigl). & (4.13) \cr
              }$$
By (3.5), $\sigma(Y_\ell) \subset 
\sigma(\overline Y_\ell, U_\ell)$ for each 
integer $\ell$.
One now has by (3.8) and (4.11) that
$$ 
\sigma(\Upsilon_\ell, \ell \geq N)
\subset \sigma({\bf B}). \eqno (4.14) $$
We need some sort of analog of (4.14) for
$\sigma(\Upsilon_\ell, \ell \leq 0)$ and
$\sigma({\bf A})$.

   For that purpose, we will need to work with the
probability measure $P_0$ (in (3.6)).   
Some more notations will be needed.
For events $A$ and $B$, the notation $A \doteq B$ will
mean that $P_0(A \triangle B) = 0$, where $\triangle$
denotes the symmetric difference.
For an event $A$ and a $\sigma$-field ${\cal B}$, the
notation $A \dot\in {\cal B}$ will mean that
$A \doteq B$ for some $B \in {\cal B}$.
For $\sigma$-fields ${\cal A}$ and ${\cal B}$, the 
notation ${\cal A} \dot\subset {\cal B}$ will mean
that $A \dot\in {\cal B}$ for every 
$A \in {\cal A}$, and the notation
${\cal A} \doteq {\cal B}$ will mean that
${\cal A} \dot\subset {\cal B}$ and
${\cal B} \dot\subset {\cal A}$.

   Refer to (3.7) and both equalities in (3.6).
For $\omega \in \{\xi_0 = 0\}$, the ``ordered pairs''
$(V_k(\omega), k \leq 0)$ determine (``measurably'')
the set of integers
$\{\xi_k(\omega), k \leq 0\}$ as well as 
$Z_j(\omega)$ ($M$ or $I$) for $j$ in that set,
and hence determine $Z_j(\omega)$ for all $j \leq 0$   
(since $Z_j(\omega) = D$ for integers
$j \leq 0$ that are not in that set).
Combining that with (3.3), one obtains that
$\sigma(V_k, k \leq 0) \doteq \sigma(Z_k, k \leq 0)$.
Now  
$\sigma(\overline X_k, k \leq 0) \subset
\sigma(X_k, Z_k, k \leq 0)$
by (3.1) and (3.2); hence
$\sigma(\overline Y_k, k \leq 0) \subset
\sigma(X_k, Z_k, k \leq 0)$ 
by (3.3) and (3.4); and hence also 
$\sigma(Y_k, k \leq 0) \subset
\sigma(X_k, Z_k, U_k, k \leq 0)$ by (3.5).
Hence by (4.12) and (3.8),
$$ \sigma(\Upsilon_k, k \leq 0) \dot\subset 
\sigma({\bf A}). \eqno (4.15) $$

   {\it Step 5.} 
On the original probability space $(\Omega, {\cal F} ,P)$
let $\Gamma$ be a random variable which (under $P$) is uniformly distributed on the interval $[0,1]$ and is
independent of $\sigma(X,Z,U)$ (recall paragraphs (P2),
(P3), and (P8)).
Now of course by the hypothesis (of Lemma 4.4), the sequence
$\beta(X,n;P) \to 0$ as $n \to \infty$
(absolute regularity).
At this point, we shall apply the ``coupling'' result of
Berbee [1, p.\ 104, Theorem 4.4.7 and p.\ 106, the Note],
which is closely related to the ``maximal coupling''
result of Goldstein [9].
As a convenient reference for Berbee's result, we cite 
[5, v2, p.\ 277, Theorem 20.7 and
pp.\ 477-478, Lemma A1651].
The latter ``lemma'' in that reference simply involves
the use of the random variable $\Gamma$ above as a
``randomizer'', and it is simply (a special case of)
the version in Dudley and Philipp [8, Lemma 2.11]
of a theorem of Skorohod [15].
Thereby there exists a sequence 
$X' := (X'_k, k \in {\bf Z})$
of random variables with the following properties
(under $P$):
$$ \eqalignno{
& {\rm the\ distributions\ of\ the\ sequences}\ X\  
{\rm and}\ X'\ 
{\rm are\ identical\ (under}\ P); & (4.16) \cr
& {\rm the}\ \sigma{\rm-fields}\  
\sigma(X')\ {\rm and}\ \sigma(X_k, k \leq 0)\
{\rm are\ independent\ (under}\ P); & (4.17)\cr
&\forall n \in {\bf N},
P(\exists k \geq n\ {\rm such\ that}\
X'_k \neq X_k) = \beta(X,n;P);\ {\rm and} & (4.18) \cr 
& \sigma(X') \subset \sigma(X,\Gamma). & (4.19) \cr
}$$
By paragraphs (P1)-(P8) and the properties of the random variable $\Gamma$,
one has that under $P$ the $\sigma$-fields 
$\sigma(\Gamma)$, $\sigma(X)$, $\sigma(Z)$, and 
$\sigma(U)$ are independent.
By (3.6) and a trivial argument, that holds under $P_0$
as well.
Further, by (3.6) and a trivial argument, 
the distribution of the random array
$(\Gamma, X, X', U)$ 
(on $[0,1] \times {\bf R}^{\bf Z}
\times {\bf R}^{\bf Z} \times [0,1]^{\bf Z}$)
is the same under $P_0$ as it is under $P$.
In particular, (4.16)-(4.18) all hold with $P$ replaced
with $P_0$ (see also (3.9)).

   Thus under $P_0$, the following statements hold:
By (4.19) the $\sigma$-fields,
$\sigma(X,X')$, $\sigma(Z)$, and $\sigma(U)$
are independent;
hence by (4.17) the $\sigma$-fields
$\sigma(X')$, $\sigma(X_k, k \leq 0)$, 
$\sigma(Z)$, and $\sigma(U)$ are independent;
and hence the $\sigma$-field $\sigma(X')$ is 
independent of the $\sigma$-field 
$\sigma(U,Z) \vee \sigma(X_k, k \leq 0)$.
\bigskip

   {\it Step 6.}
Now referring to (4.5), define analogously to (4.6)
the (``one-sided'') sequence
$X'^* := (X_1'^*, X_2'^*,X_3'^*, \dots)$ 
of random variables as follows:
$$ \forall k\geq 1,\ X_k'^* := X'_{k+T}. \eqno (4.20) $$

   Now consider an arbitrary event 
$A \subset \sigma(U,Z) \vee \sigma(X_k, k \leq 0)$
such that $P_0(A) >0$.  
For each integer $t \geq 0$ such that 
$P_0(A \cap \{T=t\}) > 0$, 
one now has (note that $\sigma(T) \subset \sigma(Z)$) that
$$ \eqalignno{ {\cal L}_0(X'^*|A \cap \{T = t\}) &=
{\cal L}_0((X'_{t+1}, X'_{t+2}, X'_{t+3}, \dots)|
A \cap \{T=t\}) \cr
&= {\cal L}_0(X'_{t+1}, X'_{t+2}, X'_{t+3}, \dots)
= {\cal L}_0(X'_1, X'_2, X'_3, \dots), & (4.21) \cr 
}$$
where ${\cal L}_0(\dots)$ 
(resp.\ ${\cal L}_0(\dots|\dots)$)
denotes the distribution (resp.\ conditional
distribution) under $P_0$.
Since the last term in (4.21) is ``constant''
(not depending on $A$ or $t$),
it follows by a simple standard calculation that
for each such event $A$,
${\cal L}_0(X'^*|A) =  
{\cal L}_0(X'_1, \allowbreak X'_2, X'_3, \dots)$,
and also (consider the case $A = \Omega$) 
${\cal L}_0(X'^*) =  
{\cal L}_0(X'_1, \allowbreak X'_2, X'_3, \dots)$.
Consequently, the sequence $X'^*$ is 
(under $P_0$) independent of the 
$\sigma$-field
$\sigma(U,Z) \vee \sigma(X_k, k \leq 0)$.         
                                                 
   Analogously to (4.13),  
define the random array $\bf B'$ as follows:
$$  
{\bf B'} := \Bigl(X'^*; (V_k, k \geq N);
              (U_k, k \geq N)\Bigl). \eqno (4.22) 
$$ 
Referring to the last sentence of the preceding 
paragraph and the third sentence after (4.19)
(which together with paragraph (P8)(i) yields
$\beta(U,N;P_0) = 0$), 
one has by (4.12), (4.22), 
Remark 3.1(B), and Lemma 2.1 that  
$$\beta(\sigma({\bf A}), \sigma({\bf B'});P_0) = \beta(V,N;P_0). \eqno (4.23) $$
Also, by (4.5), (4.13), (4.22), and (4.18)
(and the fact $\sigma(T) \subset \sigma(Z)$), 
with the sums below taken over all nonnegative 
integers $t$ such that
$P_0(T=t) > 0$, one has 
(recall the sequence $H(\, .\, )$ in the statement of
Lemma 4.4) that 
$$ \eqalignno{ 
P_0({\bf B'} &\neq {\bf B}) =
P_0(X'^* \neq X^*) =
\sum P_0(X'^*\neq X^*|T = t) 
\cdot P_0(T = t) \cr
&= \sum 
P_0\bigl((X'_{t+1}, X'_{t+2}, X'_{t+3}, \dots) \neq       
(X_{t+1}, X_{t+2}, X_{t+3}, \dots)|T=t \bigl) 
\cdot P_0(T=t) \cr
& =  \sum  
P_0\bigl((X'_{t+1}, X'_{t+2}, X'_{t+3}, \dots) \neq       
  (X_{t+1}, X_{t+2}, X_{t+3}, \dots)\bigl) \cdot P_0(T=t) \cr   
& = \sum \beta(X,t+1;P_0) \cdot P_0(T=t) \cr
&= E_0H(T). \cr   
}$$
Hence by (4.12)-(4.15), (4.23), and Lemma 2.2
(and Theorem 4.1(I), proved above) 
$$ \eqalignno{ 
\beta(\Upsilon, N; P_0) &\leq  
\beta(\sigma({\bf A}), \sigma({\bf B}); P_0) 
\leq \beta(\sigma({\bf A}), \sigma({\bf B'}); P_0) + E_0H(T)
\cr
&= \beta(V,N; P_0) + 2E_0H(T). & (4.24) \cr 
}$$

   Now by (4.5), 
$T \geq \sum_{i=1}^{N-1} I(Z_{\xi(i)} = M)$.
Also, the sequence $H(n)$ (in the statement of Lemma 4.4)
is nonincreasing as $n$ increases.
It follows that
$E_0H(T) \leq E_0H(\sum_{i=1}^{N-1} I(Z_{\xi(i)} 
\allowbreak = M))$.
Combining this with (4.24) one obtains Lemma 4.4.
\bigskip

   {\it Proof of Statement (II) in Theorem 4.1.}
Define the sequence 
$\Theta := (I(Z_{\xi(i)}= M), i \in {\bf Z})$
of random indicator functions.
Under $P_0$ this sequence is strictly stationary
by (3.7) and Theorem 4.1(I) (proved above).
By (3.6), paragraph (P3), and a trivial argument,
the sequence $\Theta$ is also nondegenerate under $P_0$.
Also, by (3.7) and the hypothesis (of Theorem 4.1(II)), 
one has that
$ \sum_{n=1}^\infty \beta(\Theta,n; P_0) < \infty$.
Also, by the hypothesis (of Theorem 4.1(II)),
the (nonincreasing) sequence 
$(H(n), n \in \{0\} \cap {\bf N})$ of
nonnegative numbers in Lemma 4.4 is summable. 
Hence for that sequence $H(\, .\, )$, by Lemma 2.3,
$\sum_{n=2}^\infty 
E_0H(\sum_{i=1}^{n-1} I(Z_{\xi(i)} 
\allowbreak = M)) < \infty$.  
Hence by (3.8), Lemma 4.4, and the hypothesis 
(of Theorem 4.1(II)), 
the conclusion of Theorem 4.1(II) holds.
That completes the proof.
\bigskip

\centerline {\bf References}
\bigskip
\def\refs{\medskip\hangindent=25pt\hangafter=1\noindent}

\refs {[1]} Berbee, H.C.P.\ (1979).  
{\it Random Walks with Stationary Increments and 
Renewal Theory\/}.
Mathematical Centre, Amsterdam.

\refs {[2]} Billingsley, P.\ (1995).  
{\it Probability and Measure\/}, 3rd edn.
Wiley, New York.

\refs {[3]} Bradley, R.C.\ (1980).
A remark on the central limit question for dependent
random variables.
{\it J.\ Appl. Probab.\/} {\bf 17}, 94-101.

\refs {[4]} Bradley, R.C.\ (1989).
A stationary, pairwise independent, absolutely
regular sequence for which the central limit theorem
fails.
{\it Z.\ Wahrsch.\ verw.\ Gebiete\/} {\bf 81}, 1-10. 

\refs {[5]} Bradley, R.C.\ (2007).
{\it Introduction to Strong Mixing Conditions\/},
Volumes 1, 2, and 3.
Kendrick Press, Heber City (Utah).

\refs {[6]} Chaudhuri, P.\ and Dasgupta, A.\ (2006).
Stationarity and mixing properties of replicating
character strings.
{\it Statistica Sinica\/} {\bf 16}, 29-43.

\refs {[7]} Davydov, Yu.A.\ (1973).
Mixing conditions for Markov chains.
{\it Theor.\ Probab.\ Appl.\/} {\bf 18}, 312-328.

\refs {[8]} Dudley, R.M.\ and Philipp, W. (1983).  
Invariance principles for sums of Banach space valued 
random elements and empirical processes. 
{\it Z.\ Wahrsch.\ verw.\ Gebiete\/} {\bf 62}, 509-552.

\refs {[9]} Goldstein, S.\ (1979).
Maximal coupling.
{\it Z.\ Wahrsch.\ verw.\ Gebiete\/} {\bf 46}, 193-204.

\refs {[10]} Ibragimov, I.A.\ and Linnik, Yu.V.\ (1971).
{\it Independent and Stationary Sequences of Random 
Variables\/}.
Wolters-Noordhoff, Groningen.

\refs {[11]} Merlev\`ede, F.\ and Peligrad, M.\ (2000).  
The functional central limit theorem under the strong 
mixing condition.  
{\it Ann.\ Probab.\/} {\bf 28}, 1336-1352.

\refs {[12]} Pinsker, M.S.\ (1964).  
{\it Information and Information Stability of Random 
Variables and Processes\/}.
Holden-Day, San Francisco.

\refs {[13]} Rio, E.\ (2000).  
{\it Th\'eorie asymptotique des
processus al\'eatoires faiblement d\'ependants\/}. 
Math\'ematiques \& Applications 31.  
Springer, Berlin.

\refs {[14]} Rosenblatt, M.\ (1956).
A central limit theorem and a strong mixing condition.
{\it Proc.\ Natl.\ Acad.\ Sci.\ USA\/} {\bf 42}, 43-47.

\refs {[15]} Skorohod, A.V.\ (1976).  
On a representation of random variables.  
{\it Theor.\ Probab.\ Appl.\/} {\bf 21}, 628-632.

\refs {[16]} Volkonskii, V.A.\ and Rozanov, Yu.A.\ (1959).
Some limit theorems for random functions I.
{\it Theor.\ Probab.\ Appl.\/} {\bf 4}, 178-197.

\refs {[17]} Waterman, M.S.\ (1995).
{\it Introduction to Computational Biology\/}.
Chapman and Hall, New York.
\bigskip

\bye